\documentclass[12pt,twoside]{amsart}
\usepackage{mathrsfs}
\usepackage{amsfonts}
\usepackage{amssymb}
\usepackage{amsmath}
\usepackage{fancyhdr}
\usepackage{graphicx}

\theoremstyle{plain}
\newtheorem{thm}{Theorem}[section]
\newtheorem{theorem}[thm]{Theorem}
\newtheorem{lemma}[thm]{Lemma}

\newtheorem{proposition}[thm]{Proposition}

\theoremstyle{definition}

\newtheorem{remark}[thm]{Remark}

\newtheorem{definition}[thm]{Definition}

\newtheorem{defn-thm}[thm]{Definition-Theorem}

\numberwithin{equation}{section}

\begin{document}

\title{Rectifiability of singular sets for geometric flows (I)--Yang-Mills flow}

\author{Jian Zhai \\ \tiny{Department of Mathematics, Zhejiang University, Hangzhou
310027, PRC} } \maketitle


\begin{abstract}
We prove that monotonicity of density and energy inequality imply
the rectifiability of the singular sets for  Yang-Mills flow.
\end{abstract}

\section{Introduction}

In this serial of papers we study the properties of singular sets
for geometric flows, such as harmonic map heat flow, Yang-Mills
flow, mean curvature flow and Ricci flow. Particularly we are
interested in the rectifiability of singular sets and the structure
of limit flows for these geometric flows. Generally, we can not
expect a global smooth solution to these flows. In many cases it has
been proved that these geometric flows may develop singularity as
$t$ increasing (see \cite{HuS}\cite{N}\cite{T}\cite{Z}). So to
understand the properties and structure of the singularity is a very
important problem of the regularity theory of theses geometric
flows. \\

 A common property of these
geometric flows is the monotonicity of density function. A
remarkable result of these papers is that the monotonicity of
density and energy inequality imply the rectifiability of the
singular sets.  Furthermore if the monotonicity formula for density
is equality, then all tangent flows of a mean convex mean curvature
flow  are the dilation of $\mathbb R^n$. This is the best possible
result for the tangent flows (see \cite{Hu}\cite{W}). The rectifiability
and structure of singular sets for harmonic map heat flow was considered in \cite{Z}.\\

It is well-known that monotonicity methods are often used in the
regularity problems for minimal surfaces, harmonic maps and
Navier-Stokes equations, as well as other nonlinear partial
differential equations.  The monotonicity of density for Yang-Mills
flow in static case was first proved by \cite{P}. The monotonicity
inequality for a parabolic version of density function, and
$\epsilon$-regularity of a solution to Yang-Mills flow were proved
in \cite{CS}\cite{H}\cite{N}\cite{T}.
  Hamilton \cite{H} proved monotonicity formulas on manifolds for smooth parabolic flows,
  including  Yang-Mills flow, mean curvature flow and
  harmonic heat flow.\\

  In this paper we only consider the Yang-Mills flow. The rectifiability of first time singular set for Yang-Mills
  flow were proved in \cite{T}. Moser \cite{Mo} considered the rectifiability
  of singular sets for a class of static solutions including static solutions of Yang-Mills
  equations.

{\large\bf Acknowledgments.} \quad  This work was done during the
period of visiting University of California at Berkeley. I would
like to thank Professor Evans for his hospitality and enlightening
discussions. This work is supported by NSFC No.10571157 and China
Scholarship Council.

\section{Yang-Mills flow}

Yang-Mills equations as a field theory were originally arisen in
physics. Here we consider the Yang-Mills flow--a parabolic version
of Yang-Mills equations. The Yang-Mills flow was first proposed by
Atiyah-Bott\cite{AB}, and it played an important role in Donaldson's
four
dimensional topology theory \cite{D}.\\

We consider the Yang-Mills flow:
\begin{equation} \label{1.1}
\begin{aligned}
\partial_\tau A=-\triangledown^* F,\quad\text{in}\quad \mathbb R^m\times
[0,T)\\
A(x,0)=A_0(x),\quad\text{in}\quad \mathbb R^m
\end{aligned}
\end{equation}
where  $A(x,\tau)=(A_1(x,\tau), A_2(x,\tau),\dots, A_m(x,\tau))$ is
a gauge potential, and for a fixed compact Lie group $G$,  $A_\mu$
$(\mu=1,2,\dots, m)$ take values in Lie algebra $Lie(G)$. The readers may see \cite{A}\cite{T}
for the explanation of the notations. However it is not necessary to read \cite{A}\cite{T} for understanding
the proofs of this paper. In fact, it is enough to consider following special situation.   \\

 For simplicity, in this paper we may assume
$$
G=U(n):=\{\text{$b$ is an $n\times n$ complex matrix with}:\quad
bb^*=I \},
$$
where $b^*$ denotes the conjugate transpose of $b$, and $I$ is the
unit matrix. The Lie algebra $Lie(U(n))$ can be expressed as
$$
Lie(U(n))=\{\text{$B$ is an $m\times m$ complex matrix with}:\quad
B+B^*=0\}.
$$\\

To define an inner product, we need the following Killing form on
$Lie(U(n))$ defined by
\begin{equation}\label{1.01}
<B,C>=(-1)\text{trace}(BC),\quad\forall B,\,\,\, C\in Lie(U(n)).
\end{equation}
Here trace($BC$) means the trace of the matrix $BC$. Furthermore, if
$$
B(x)=(B_1(x), B_2(x),\dots,B_k(x)\quad\text{and}\quad C(x)=(C_1(x),
C_2(x),\dots,C_k(x)),$$
 and $B_\mu(x)$ and $C_\mu(x)$
 ($\mu=1,2,\dots,k$) are $Lie(U(n))$-valued functions, we define the
inner product
\begin{equation}\label{1.02}
(B,C)=\int_{\mathbb R^m}\sum_{1\leq\mu\leq
 k}<B_\mu(x),C_\mu(x)>dx=\int_{\mathbb
R^m}\sum_{1\leq \mu\leq k}(-1)\text{trace}(B_\mu(x)C_\mu(x))dx.
\end{equation}\\

 Associated with the
potential $A$, the operator $\triangledown$ and the gauge field $F$
are defined by
$$
\triangledown_\mu :=\frac{\partial}{\partial x_\mu}+A_\mu
$$
and
$$
F_{\mu\nu}=[\triangledown_\mu, \triangledown_\nu]=\frac{\partial
A_\nu}{\partial x_\mu}-\frac{\partial A_\mu}{\partial
x_\nu}+[A_\mu,A_\nu]
$$
respectively. The operator $\triangledown^*$ in (\ref{1.1}) is the adjoint operator
of $\triangledown$ with respect to the inner product (\ref{1.02}).   \\

The Yang-Mills flow (\ref{1.1}) can be regarded as a gradient flow
of an action. The invariant action of the Yang-Mills flow is defined
by
$$
YM(A)(\tau)=\int_{\Bbb R^m}|F(x,\tau)|^2dx,\quad
|F(x,\tau)|^2:=\sum_{1\leq\mu,\nu\leq
m}(-1)\text{trace}(F_{\mu\nu}(x,\tau)F_{\mu\nu}(x,\tau)).
$$\\

A parabolic version of density function to (\ref{1.1}) is defined by
\cite{CS}\cite{H}\cite{N}
$$
\theta^\rho(z,\tau)=\rho^{4-m}\int_{\mathbb
R^m}\exp\{-\frac{|z-y|^2}{4\rho^2}\}|F(y,\tau-\rho^2)|^2dy,\quad\forall\rho>0.
$$
The monotonicity of density function for static Yang-Mills flow was
first proved by \cite{P}. For the parabolic version of density
function, the monotonicity inequality
\begin{equation}
\theta^\rho(z,\tau)\leq C_0\theta^{\rho'}(z,\tau)+C_1,\quad
0<\forall \rho\leq\forall\rho',
\end{equation}
 and $\epsilon$-regularity of
a smooth solution to (\ref{1.1}) were proved in
\cite{CS}\cite{H}\cite{N}.\\

The $\epsilon$-regularity means that there is a constant
$\epsilon>0$ such that $\theta^\rho(z,\tau)<\epsilon$ for some
$\rho>0$ implies that
 $(z,\tau)$ is a regular point of $A$. So we can define the singular
set of the Yang-Mills flow as the
complement of the set of the regular points.\\

\begin{definition}\label{d1.1}For $\tau>0$, the singular set of the
Yang-Mills flow (\ref{1.1}) is defined by
\begin{equation}\label{1.05}
S(\tau):=\{z\in\mathbb R^m:\quad
\liminf_\rho\theta^\rho(z,\tau)\geq\epsilon\}.
\end{equation}
\end{definition}

\vspace{1cm}

Notice that from \cite{N}\cite{HT} etc., the singular set $S(\tau)$
is not empty for some time $\tau>0$, and the $(m-4)$-dimension
Hausdorff measure of $S(\tau)$ is finite. In this paper, we shall
prove the rectifiability of $S(\tau)$.

\begin{theorem} \label{th1} If $A(x,\tau)$ is a weak solution of
(\ref{1.1}) and its density function has monotonicity, then the
singular set $S(\tau)$ of $A(x,\tau)$ at time $\tau>0$ is ($m-4$)-
rectifiable.
\end{theorem}

\vspace{1cm}

\begin{remark}\label{rm1} (1) The rectifiability for the singular set of steady state Yang-Mills equations is proved by
\cite{Mo}, and for the first time singular set of Yang-Mills flow is
proved by \cite{T}\cite{HT}.\\
(2)If $A(x,\tau)$ is a weak solution of (\ref{1.1}), then as in
\cite{Z}\cite{Z2}, the motion of the singular set $S(\tau)$ can be
described by the limit equation of
\begin{equation}\label{1.2}
\partial_\tau\theta^\rho(z,\tau)=\Delta_z\theta^\rho(z,\tau)+\frac1{\rho^2}(\theta^\rho(z,\tau)-\frac{\rho}2\partial_\rho\theta^\rho(z,\tau)),\quad \forall (z,\tau)\in\mathbb R^m\times\mathbb R^+.
\end{equation}
 The further results on the motion of the singular set
may be obtained by considering the properties of solutions to the
Yang-Mills flow (see \cite{Z}). \\
(3) Similar results of this paper can be proved on Riemannian
manifolds instead of $\mathbb R^m$.
\end{remark}

\vspace{1cm}

\section{Energy inequality, monotonicity and local continuity of density for Yang-Mills flows }

In this section, we shall introduce several basic properties of the
Yang-Mills flow, including the energy inequality (Lemma \ref{l2.1}),
the monotonicity of the density function (Lemma \ref{l2.2}) and
local continuity of the density function (Lemma \ref{l2.3}). The
first and second lemma are known fact, and their proofs can be found
in the references. Lemma \ref{l2.3} was first proved in \cite{Z},
and for readers' convenience, we give its proof here.

\begin{lemma}\label{l2.1}(energy inequality\cite{D}\cite{CS}\cite{N}) Suppose $A$
is a smooth solution of (\ref{1.1}). Then for all $\tau\in(0,T)$
$$
YM(A)(\tau)+2\int_0^\tau\int_{\mathbb R^m}|\partial_t
A(y,t)|^2dydt=YM(A_0).
$$
Particularly we have the energy inequality
\begin{equation}\label{2.1}
YM(A)(\tau)\leq YM(A_0).
\end{equation}

\end{lemma}

\vspace{1cm}

\begin{lemma}\label{l2.2}(monotonicity \cite{CS}\cite{H}\cite{N}) Suppose $A$ is a smooth
solution of (\ref{1.1}). Then there is a constant $C>0$ such that
for all $\rho$, $\rho'$: $0<\rho\leq\rho'$, we have
\begin{equation}\label{2.2}
\theta^\rho(z,\tau)\leq
C_0(\rho,\rho')\theta^{\rho'}(z,\tau)+C_1(\rho,\rho'),\quad
\text{for}\quad z\in \mathbb R^m,\quad \tau>0,
\end{equation}
where
$$
C_0(\rho,\rho')=C\exp\{C(\rho'-\rho)\},\quad
C_1(\rho,\rho')=C(\rho'^2-\rho^2)YM(A_0).
$$
\end{lemma}

\vspace{1cm}

\begin{lemma}\label{l2.3}(local continuity) Suppose $A$ is a weak solution of (\ref{1.1})
which satisfies the energy inequality (\ref{2.1}) and monotonicity
(\ref{2.2}).  For all $R>1$, for all
$\alpha=(\alpha_1,\alpha_2,\dots,\alpha_m)\in \mathbb N^m$, there is
a constant $C$ which is independent of  $\rho>0$ and $x,\tilde x\in
B_R$, such that
\begin{equation}\label{2.10}
 |D^\alpha_{x}\theta^\rho(\bar z+\rho
x,\tau)-D^\alpha_{x}\theta^\rho(\bar z+\rho\tilde x,\tau)|\leqq
C|x-\tilde x|.
\end{equation}

\end{lemma}

$Proof.$ Note that
\begin{eqnarray*}
&&\theta^\rho(\bar z+\rho x,\tau)-\theta^\rho(\bar z+\rho\tilde x,\tau)\\
&&=\frac12\int_{\mathbb R^m}\{e^{-\frac{|x-\bar x|^2}4}-e^{-\frac{|\tilde x-\bar x|^2}4}\}|F(\bar z+\rho\bar x,\tau-\rho^2)|^2  d\bar x\\
&&=\frac12\int_{|x-\tilde x||\bar x-x|\leqq r}e^{-\frac{|x-\bar x|^2}4}\{1-e^{\frac{|x-\bar x|^2-|\tilde x-\bar x|^2}4}\} |F(\bar z+\rho\bar x,\tau-\rho^2)|^2 d\bar x\\
&&+\frac12\int_{|x-\tilde x||\bar x-x|\geqq r}e^{-\frac{|x-\bar x|^2}4}|F(\bar z+\rho\bar x,\tau-\rho^2)|^2d\bar x\\
&&-\frac12\int_{|x-\tilde x||\bar x-x|\geqq r}e^{-\frac{|\tilde
x-\bar x|^2}4}|F(\bar z+\rho\bar x,\tau-\rho^2)|^2d\bar x.
\end{eqnarray*}

Case 1. $|x-\tilde x||\bar x-x|\leqq r$ and $|x-\tilde x|\leqq\frac
r{2R}$.

We have
\begin{equation}\label{2.11}
\begin{aligned}
|1-\exp\{\frac{|x-\bar x|^2-|\tilde x-\bar x|^2}4\}|\exp\{-\frac{|x-\bar x|^2}8\}\\
=|1-\exp\{\frac{(x-\tilde x)\cdot(x+\tilde x)+2\bar x\cdot(\tilde x-x)}4\}|\exp\{-\frac{|x-\bar x|^2}8\}\\
\leq \frac{|x-\tilde x|(|x+\tilde x|+2|\bar x|)}2\exp\{-\frac{|x-\bar x|^2}8\}\\
\leq \frac{|x-\tilde x|(|x+\tilde x|+2|x|+2|\bar x-x|)}2\exp\{-\frac{|x-\bar x|^2}8\}\\
\leq C|x-\tilde x|
\end{aligned}
\end{equation}
where we used the estimates
$$
|x-\tilde x|(|x+\tilde x|+2|\bar x|)\leqq |x-\tilde x||x+\tilde
x|+2|x-\tilde x||\bar x-x|+2|x-\tilde x||x|\leqq 4r
$$
and
$$
|1-\exp\{-h\}|\leqq 2h,\quad\text{for}\quad 0\leqq h\leqq r<1,
$$
because $r$ is small and
$$
|\bar x-x|\leq\frac{r}{|x-\tilde x|}.
$$

Case 2. $|x-\tilde x||\bar x-x|\geqq r$ and $|x-\tilde x|\leqq\frac
r{2R}$.

We have
\begin{equation}\label{2.12}
\exp\{-\frac{|x-\bar x|^2}8\}\leqq\frac {C}{|x-\bar x|^2}\leqq
\frac{C|x-\tilde x|^2}{r^2},
\end{equation}
 and
\begin{equation}\label{2.13}
\exp\{-\frac{|\tilde x-\bar x|^2}8\}\leqq\frac{C}{|\tilde x-\bar
x|^2}\leqq\frac{C|x-\tilde x|^2}{(r-r^2/4)^2}
\end{equation}
 where we have used
\begin{eqnarray*}
|\bar x-\tilde x||x-\tilde x|&&\geqq (|\bar x-x|-|x-\tilde x|)|x-\tilde x|\\
&&\geqq |\bar x-x||x-\tilde x|-|x-\tilde x|^2\\
&&\geqq r-\frac{r^2}4.
\end{eqnarray*}

Note that
\begin{equation}\label{2.14}
\int_{\Bbb R^m}e^{-\frac{|x-\bar x|^2}8}|F(\bar z+\rho\bar
x,\tau-\rho^2)|^2d\bar x\leqq C(YM(A_0)+1)
\end{equation}
because if we take
$$
\bar x:=\frac{(\bar x-x)}{\sqrt{2}},\,\,\bar\rho:=\sqrt{2}\rho
$$
then
\begin{eqnarray*}
&&\int_{\Bbb R^m}e^{-\frac{|x-\bar x|^2}8}|F(\bar z+\rho\bar x,\tau-\rho^2)|^2d\bar x\\
&&=2^{\frac{m-4}2}\int_{\Bbb R^m}e^{-\frac{|\bar
x|^2}4}|F(\bar z+\rho x+\bar \rho\bar x,\tau+\rho^2-\bar\rho^2)|^2d\bar x \quad(\text{by $x-\bar x\to \sqrt{2}\bar x$})\\
&&=2^{\frac{m-4}2}\theta^{\bar\rho}(\bar z+\rho x,\tau+\rho^2)\\
&&\leq 2^{\frac{m-4}2}\{C_0(\bar\rho,\sqrt{\tau})\theta^{\bar\rho}(\bar z+\rho x,\tau+\rho^2)+C_1(\bar\rho,\sqrt{\tau})\},\quad(\text{by monotonicity})\\
&&\leqq C(YM(A_0)+1)\quad(\text{by the energy inequality}).
\end{eqnarray*}

Fix $r$ (small). By (\ref{2.11})-(\ref{2.14}) we obtain (\ref{2.10})
with $|\alpha|=0$ in the case of $x,\tilde x\in B_R$ and $|x-\tilde
x|\leqq r/2R$. It is clear for $|x-\tilde x|>r/2R$. Similarly, we
can prove (\ref{2.10}) for $|\alpha|>0$.\qed

\vspace{1cm}

\section{Rectifiability of singular sets for Yang-Mills flow}

For $k=1,2,\dots,m-1$, let $G(m,k)$ denote the set of all
$k$-dimensional subspaces passing the origin of $\mathbb R^m$.

\begin{lemma}\label{l4.2} For $\bar z\in\Bbb R^m$, we have
\begin{equation}\label{4.1}
\liminf_\rho\int_{U}\rho^{-4}\theta^{\rho}(\bar z+Y,\tau)dY<\infty
\end{equation}
 for $\gamma_{m,4}$-almost all $U\in G(m,4)$. Here
$\gamma_{m,4}$ denotes the measure on $G(m,4)$ from the
decomposition of the Lebesgue measure $dz=d\gamma_{m,4}dY$ on
$\mathbb R^m$ (see \cite{M} chapter 3).
\end{lemma}

$Proof.$ From Fubini's Theorem and Lebesgue-Fatou lemma,
\begin{eqnarray*}
\liminf_\rho\rho^{-4}\int_{\Bbb R^m}\theta^{\rho}(z,\tau)dz\\
=\liminf_\rho\int_{G(m,4)}d\gamma_{m,4}(U)\int_{U}\rho^{-4}\theta^{\rho}(\bar z+Y,\tau)dY\\
\geqq\int_{G(m,4)}d\gamma_{m,4}(U)\liminf_\rho\int_{U}\rho^{-4}\theta^{\rho}(\bar
z+Y,\tau)dY.
\end{eqnarray*}

Note that from the definition of $\theta^\rho(z,\tau)$, for all
$\tau>\rho^2$
\begin{eqnarray*}
\rho^{-4}\int_{\mathbb
R^m}\theta^{\rho}(z,\tau)dz=\rho^{-m}\int_{\mathbb
R^m}dz\int_{\mathbb
R^m}\exp\{-\frac{|z-y|^2}{4\rho^2}\}|F(y,\tau-\rho^2)|^2dy\\
= \int_{\mathbb R^m}|F(y,\tau-\rho^2)|^2dy\int_{\mathbb
R^m}\exp\{-\frac{|x|^2}4\}dx\leq CYM(A_0),
\end{eqnarray*}
by the energy inequality (\ref{2.1}).

Then for $\gamma_{m,4}$- almost all $U\in G(m,4)$,
$$
\liminf_\rho\int_{U}\rho^{-4}\theta^{\rho}(\bar z+Y,\tau)dY
$$
is finite.\qed

\vspace{1cm}

To prove the rectifiability of $S(\tau)$, we need the notation of a
cone $X(\bar z,V,s)$ around a hyperplane $V$ at $\bar z\in S(\tau)$
with cone angle $\cos^{-1} s$. For $\bar z\in\Bbb R^m$, $0<s<1$,
$0<r<\infty$ and $V\in G(m,k)$, we define the cone  as
$$
X(\bar z,V,s):=\{z\in\Bbb R^m:\quad \text{dist}(z-\bar z,V)<s|z-\bar
z|\}
$$
and define the restriction of the cone $X(\bar z,V,s)$ in an
$m$-dimensional ball $B^m(\bar z,r)$ as
$$
X(\bar z,r,V,s):=X(\bar z,V,s)\cap B(\bar z,r).
$$
Notice that $\bar z\not\in X(\bar z,V,s)$.

\vspace{1cm}

\begin{lemma}\label{l4.3} Let $\bar z\in S(\tau)$. If for some $W\in
G(m,1)$ and $\delta>0$,
\begin{equation}\label{4.2} \cap_{0<r<\delta}(S(\tau)\cap
X(\bar z,r,W,s))\not=\emptyset,\quad 0<\forall s<1,
\end{equation}
 then there is a sequence $\{\rho_i\}_i$ and a
constant $C=C(\epsilon)>0$ depending on $\epsilon$ (see \ref{1.05})
) such that
\begin{equation}\label{4.3}
\inf_{t\geq 0}\lim_{i}\theta^{\rho_i}(\bar z+\rho_i
t\omega,\tau)\geq C(\epsilon)
\end{equation}
 where $\omega\in
S^{m-1}$ is the direction of the line $W =\{t\omega,\quad \forall
t\in\Bbb R\}$.
\end{lemma}

$Proof.$ Fix $\bar z\in S(\tau)$.  Note that  for any $t> 0$ and
$h>0$, if we take $s$ such that
 $$
 s\leq \frac{h/2}{\sqrt{t^2+(h/2)^2}},
 $$
 then
 $$
 \cup_{\rho_0\geq \rho>0}B^m(\bar z+\rho t\omega,\rho h/2)\supset
 X^+(\bar z,r,W,s),
 $$
here $+$ means the positive part of the cone $X$. Note that if we
take $t< 0$, then we can cover the negative part of the cone $X$.

 So
from (\ref{4.2}) there is a sequence $\{\rho_i\}_i$ such that for
each $i\in\Bbb N$,
$$
S(\tau)\cap B^m(\bar z+\rho_it\omega,\rho_ih/2)\not=\emptyset.
$$
Take $\bar z_{\rho_i}$ belonging to this joint set. Then  $\bar
z_{\rho_i}\in S(\tau)$ satisfies
\begin{equation}\label{4.4}
\bar z_{\rho_i}\to\bar z,\quad\text{and}\quad B^m(\bar
z+\rho_it\omega,\rho_i h)\supset B^m(\bar z_{\rho_i},\rho_i h/2).
\end{equation}

We may suppose that for $\{\rho_i\}_i$  the limit
$$
\lim_i\theta^{\rho_i}(\bar z+\rho_ix,\tau), \quad\forall x\in\Bbb
R^m
$$
exists (if not, we can take a subsequence).

So we have
\begin{equation}\label{4.5}
\begin{aligned}
\int_{B^m(t\omega,h)}\lim_i\theta^{\rho_i}(\bar z+\rho_ix,\tau)dx \\
=\lim_{i}\rho_i^{-m}\int_{B^m(\bar z+\rho_it\omega,\rho_i h)}\theta^{\rho_i}(z,\tau)dz \quad(\text{by $z=\bar z+\rho_i x$})\\
\geqq\lim_i\rho_i^{-m}\int_{B^m(\bar z_{\rho_i},\rho_i h/2)}\theta^{\rho_i}(z,\tau)dz \quad(\text{by (\ref{4.4})})\\
\geqq\inf_{ z\in S(\tau)}\int_{B^m(0,h/2)}\liminf_\rho\theta^{\rho}(z+\rho x,\tau)dx \quad(\text{by $z\to \bar z_{\rho_i}+\rho_i x$})\\
\geq C(\epsilon)h^m, \quad\text{(by (\ref{1.05}), where
$C(\epsilon)>0$ only depends on $\epsilon$)}.
\end{aligned}
\end{equation}
 Noting the continuity of
$\theta^\rho(\bar z+\rho x,\tau)$ about $x$ (see Lemma \ref{l2.3}),
the left side of (\ref{4.5}) is smaller than
$$
Ch^m\lim_i\theta^{\rho_i}(\bar z+\rho_it\omega,\tau)+C(t)h^{m+1}.
$$
If we fix $h$ such that
$$
C(\epsilon)-C(t)h\geq C(\epsilon)/2
$$
then we get (\ref{4.3}).\qed

\vspace{1cm}

\begin{proposition}\label{pro4.1} If $A(x,\tau)$ is a weak solution of
(\ref{1.1}) which satisfies the energy inequality and its density
function satisfies the monotonicity inequality, then the  singular
set $S(\tau)$ of $A(x,\tau)$ at time $\tau>0$ is ($m-4$)-
rectifiable.
\end{proposition}

$Proof$ If $S(\tau)$ is not $(m-4)$-rectifiable, then there is a
purely $(m-4)$-unrectifiable subset $S$ of $S(\tau)$ such that for
any $U\in G(m,4)$, for any $\delta>0$
\begin{equation}\label{4.6}
\limsup_{s\to0}\sup_{0<r<\delta}(rs)^{4-m}\mathcal{H}^{m-4}(S\cap
X(\bar z,r,U,s))>0
\end{equation}
 for $\mathcal{H}^{m-4}$ almost all $\bar z\in S$ (see \cite{M} Corollary
 15.15).\\

 We select special sequences of $\delta$, $s$ and $r$ by the following three
 steps.

(1) Take a sequence $\{\delta_i\}_{i\in\Bbb N}$ such that
$\delta_i\downarrow 0$ as $i\uparrow \infty$;

(2) For each $\delta_i$, take a sequence $\{s^i_j\}_{j\in\Bbb N}$
such that
$$
\text{for each $i$}\quad s^i_j\downarrow 0,\quad\text{as well
as}\quad s^j_j\downarrow 0,\quad\text{as}\quad j\uparrow\infty
$$
and
$$
\sup_{0<r<\delta_i}(rs^i_j)^{4-m}\mathcal{H}^{m-4}(S\cap X(\bar
z,r,U,s^i_j))\geqq C(i)(>0,\,\,\text{independent of $j$});
$$

(3) Furthermore, for each $i$ there is $\{r^i_j\}_{j\in\Bbb N}$:
$0<r^i_j<\delta_i$, such that
$$
 r^j_j\downarrow 0,\quad\text{as}\quad j\uparrow\infty,
$$
and for each $j$
$$
\mathcal{H}^{m-4}(S\cap X(\bar z,r^i_j,U,s^i_j))\geqq
\frac{C(i)}2(r^i_js^i_j)^{m-4}.
$$

It implies that in case of $i=j=k$ there are $\{z(k)\in S\cap X(\bar
z,r^k_{k},U,s^k_{k})\}_{k\in\Bbb N}$  such that
$$
z(k)\to \bar z, \quad\text{as}\quad k\to\infty,  \quad\text{because
$r^k_k\downarrow 0$.}
$$
Let
$$
(\omega_{1}(k),\omega_{2}(k)):=\frac{z(k)-\bar z}{|z(k)-\bar z|},
$$
where $\omega_{2}(k)\in S^3\subset U$ and $\omega_{1}(k)\to 0$ as
$k\to\infty$, because $s^k_k\downarrow 0$. By the compactness of
$S^3$, there is a subsequence (denoted also by $\{z(k)\}$) such that
$$
\omega_{2}(k)\to \omega\in S^3\subset U.
$$
Let $W=\{t\omega,\quad \forall t\in\Bbb R\}$. Then $W\subset U$ and
for any fixed $s\in (0,1)$, there is $k_0>0$ such that for $k\geq
k_0$,
$$
s^k_k\leq s,  \quad z(k)\in S\cap X(\bar z,r^k_{k},W,s),\quad
z(k)\to\bar z\quad(k\to\infty). $$

Thus, there is
$$
W\subset U,\quad W\in G(m,1)
$$
such that (\ref{4.2}) is satisfied. So by Lemma 4.2 we have
(\ref{4.3}) at least in one direction of $W$. For any $U$  in
$G(m,4)$, from the continuity of $\lim_i\theta^{\rho_i}(\bar
z+\rho_ix,\tau)$ on $x$, we get a $m$-dimension set of $x$ with
infinity $m$-dim Lebesgue measure such that on this
 set
$$
\lim_i\theta^{\rho_i}(\bar z+\rho_ix,\tau)\geqq C(>0),
$$
which is contradiction with Lemma \ref{l4.2} where $\rho$ is
replaced by $\rho_i$.\qed

\end{document}